%words123.tex: 
%%a Plain TeX file by Nathaniel Shar and Doron Zeilberger (x pages)
%begin macros

\baselineskip=14pt
\parskip=10pt

\font\eightrm=cmr8 

\magnification=\magstephalf
\def\W{{\cal W}}

\def\1{{\overline{1}}}
\def\2{{\overline{2}}}
\parindent=0pt
\overfullrule=0in

\def\frac#1#2{{#1 \over #2}}
%\headline={\rm  \ifodd\pageno  \RightHead  \else  \LeftHead  \fi}
%\def\RightHead{\centerline{
%Title
%}}
%\def\LeftHead{ \centerline{Doron Zeilberger}}
%end macros
\centerline
{\bf The  (Ordinary) Generating Functions Enumerating 123-Avoiding Words with r occurrences  }
\centerline
{\bf   of each of 1,2, ..., n are  Always Algebraic}
\bigskip
\centerline
{\it Nathaniel SHAR and Doron ZEILBERGER}

{\bf Abstract}: The set of $123$-avoiding {\it permutations} (alias words in $\{1, ..., n\}$ with exactly $1$ occurrence
of each letter) is famously enumerated by the ubiquitous Catalan numbers, whose generating function $C(x)$
famously satisfies the algebraic equation $C(x)=1+xC(x)^2$. Recently, Bill Chen, Alvin Dai, and Robin Zhou found (and very elegantly proved) an
algebraic equation satisfied by the generating function enumerating $123$-avoiding words with
{\bf two} occurrences of each of $\{1, \dots , n\}$. 
Inspired by the Chen-Dai-Zhou result,  we present an algorithm for finding such an algebraic equation for the
ordinary generating function enumerating $123$-avoiding words with exactly $r$ occurrences of
each of $\{1, \dots, n \}$ for {\bf any} positive integer $r$, thereby proving that they are {\it algebraic},
and not merely $D$-finite (a fact that is promised by WZ theory).
Our algorithm consists of presenting an {\it algebraic enumeration scheme}, combined with the Buchberger algorithm.

{\bf Introduction}

Recall that a word $w=w_1 \dots w_n$ in an ordered alphabet contains a {\it pattern} $\sigma$ (a certain permutation of $\{1, ..., k\}$)
if there exist
$$
1 \leq i_1 < i_2 < \dots < i_k \leq n \quad 
$$
such that the subword $w_{i_1} \dots w_{i_k}$ is {\it order isomorphic} to $\sigma$; in other words
$w_{i_1}, \dots, w_{i_k}$ are distinct, and if you replace the smallest entry by $1$, the second smallest entry by $2$, etc.,
you would get $\sigma$.

For example, the word {\it mathisfun} contains the pattern $132$, since (inter alia) the subword $hsn$ is order-isomorphic to $132$
(under the usual lexicographic order).

A word $w$ avoids the pattern $\sigma$ if it does not contain it. One is interested in
enumerating words, of a given length and given alphabet-size, avoiding one or more patterns. 

In a remarkable PhD thesis, under the guidance of guru Herbert S. Wilf, Alexander Burstein ([Bu]) initiated
the study of {\it forbidden patterns} (alias {\it Wilf classes}) in {\it words}, extending the
very active and fruitful research on {\it forbidden patterns} in {\it permutations}  initiated by
Donald Knuth, Rodica Simion, Richard Stanley, Herbert Wilf,  and others. For the current
{\it state of the art} of the latter, see [Wiki]. Burstein's pioneering thesis was extended by
quite a few people, and the current knowledge is described in the lucid and insightful research monographs [HM] and [Ki]. 
A systematic approach for computer-assisted enumeration of words avoiding a given set of patterns,
extending the work of Zeilberger and Vatter for permutations (see [Z4] and its references),
was initiated by
Lara Pudwell ([P]). Some of the recent work (e.g. [GGHP]) is phrased in the equivalent language
of  {\it ordered set partitions}. This equivalence is cleverly used in Anisse Kasraoui's ([Ka])
recent article.

Most of this work concerns the set of {\it all}  words avoiding a pattern.  In a very interesting recent paper
[GGHP], the authors consider (in the equivalent language of ordered set partitions), among other problems,
the problem of enumerating $123$-avoiding words of length $2n$ where each of the $n$ letters $\{1, 2, \dots, n\}$
occurs exactly twice, and conjectured a certain second-order linear recurrence equation with polynomial coefficients.
They apparently did not realize that, in their case, it was possible to justify it by a (fully rigorous, or
at least rigorizable) {\it hand-waving} argument. By general `holonomic nonsense' ([Z1]) it is known beforehand that
there is {\it some} such linear recurrence, and it is possible to bound the order, thereby justifying,
{\it a posteriori}, the guessed recurrence, provided that it is checked for sufficiently many initial values.
A more direct proof was given by Chen, Dai, and Zhou ([CDZ]), who proved the {\bf stronger} statement that
the generating function is {\it algebraic}, and even found the defining equation explicitly:
$1- \left( 2\,x + 1 \right) {F}^{2}+x \left( x+4 \right) {F}^{4}=0 \, $ .

Using Comtet's algorithm ([Co], see also [S]) for deducing, out of the algebraic equation, a linear 
differential equation for the generating function, and
hence a linear recurrence for the sequence itself, Chen, Dai and Zhou proved the [GGHP] conjecture directly.

In the present article  we will generalize this and prove that, for {\bf every} positive integer $r$, 
the ordinary generating function enumerating
$123$-avoiding words of length $rn$ where each of the $n$ letters of $\{1, 2, \dots, n\}$ occurs exactly $r$ times, is
algebraic, and present an algorithm for finding the defining equation. Alas, since at the end it uses
the memory-heavy, and exponential time, Buchberger's algorithm for finding
Gr\"obner bases, our computer (running Maple) only agreed to explicitly find the next-in-line, the analogous
equation for $r=3$:
$$
 \left( 4\,x+1 \right) ^{2}+ \left( 64\,{x}^{2}+48\,x-1 \right) {F}^{2}-2\,x \left( 128\,{x}^{2}+108\,x+27 \right) {F}^{4}-16\,{x}^{2} \left( 32\,x+27
 \right) {F}^{6}+{x}^{2} \left( 32\,x+27 \right) ^{2}{F}^{8}=0 \quad .
$$
This took less than a second, but the case $r=4$ took about an hour. Here is the minimal algebraic equation satisfied by 
the generating function, let's call it $F$, whose coefficient of $x^n$ is the number of  $123$-avoiding words with $4n$ letters
with $4$ occurrences of each $i$ ($1 \leq i \leq n$):
$$
{x}^{3}  \, \left( 5\,x-256 \right) ^{4} \left( 4\,x+1 \right) ^{4}{F}^{16}
$$
$$
+  4 \,{x}^{3} \, \left( 85\,x+58 \right)  \left( 5\,x-256 \right) ^{3} \left( 4\,x+1 \right) ^{3}{F}^{14}
$$
$$
+2\, {x}^{2} \, \left( 200\,{x}^{4}+11845\,{x}^{3}+8658\,{x}^{2}+6503\,x+256 \right)  \left( 5\,x-256 \right) ^{2} \left( 4\,x+1 \right) ^{2}{F}^{12}
$$
$$
+4\, {x}^{2} \, \left( 5\,x-256 \right)  \left( 4\,x+1 \right)  \left( 25500\,{x}^{5}-977800\,{x}^{4}+15739435\,{x}^{3}+9911721\,{x}^{2}+2082455\,x+138496 \right) 
{F}^{10}
$$
$$
+ x \, ( 60000\,{x}^{8}+2772000\,{x}^{7}-471787725\,{x}^{6}+11351360680\,{x}^{5}+15348867846\,{x}^{4}
$$
$$
+7091445146\,{x}^{3}+1387805641\,{x}^{2}
+96468480\,x-458752 ) {F}^{8}
$$
$$
+4\, x \, \left( 127500\,{x}^{7}-6439500\,{x}^{6}+28100475\,{x}^{5}+187145995\,{x}^{4}+58215739\,{x}^{3}-5955159\,{x}^{2}-2743199\,x-108800 \right) {F}^{6}
$$
$$
+ ( 10000\,{x}^{8}+628250\,{x}^{7}-57924600\,{x}^{6}+1098116930\,{x}^{5}+827342646\,{x}^{4}
$$
$$
+223797652\,{x}^{3}+
24970546\,{x}^{2}+842512\,x+1024 ) {F}^{4}
$$
$$
+ \left( 42500\,{x}^{7}-1521500\,{x}^{6}-6516800\,{x}^{5}-7480160\,{x}^{4}-276672\,{x}^{3}+461716\,{x}^{2}+
49271\,x-1024 \right) {F}^{2}
$$
$$
+ x \, \left( x+1 \right) ^{2} \left( 25\,{x}^{2}+65\,x+11 \right) ^{2} \, = 0 \, .
$$
We didn't even try the case $r=5$.

However, since we know, once again (now even without using Zeilberger's holonomic theory) that the generating function
is $D$-finite, since it has the stronger property of being algebraic, it justifies {\bf rigorously} guessing
a linear recurrence equation with polynomial coefficients, which enables one to compute, in {\it linear time},
any term of the enumerating sequence. We succeeded, using our algorithm, 
to be described below (which in particular enables a very fast enumeration of many terms of the enumerating sequences),
in discovering such recurrences for $1 \leq r \leq 5$, and using [Z2] we (or rather our beloved servant, Shalosh B. Ekhad, running
Maple) found precise asymptotics for these cases. This enables us to formulate the following intriguing conjecture,
and the second-named author (DZ) is pledging a \$100 donation to the OEIS foundation, in honor of the first prover.

{\bf Conjecture}: Let $w_r(n)$ be the number of $123$-avoiding words of length $rn$ with $r$ occurrences of
each of $\{1, \dots, n\}$. Then
$$
\lim_{n \rightarrow \infty} \,\, \frac{w_r(n)}{w_r(n-1)} = (r+1)\,2^r \quad \quad .
$$
More strongly, $w_r(n)$ is asymptotically $C_r \, \cdot \, ((r+1)2^r)^n \cdot n^{-3/2}$, where $C_r$ is a `nice' constant
(probably $\frac{1}{\sqrt{\pi}}$ times the square-root of a rational number that depends `nicely' on $r$).

Using the Maple package {\tt Words123} accompanying this article, we proved it for $r \leq 5$
(but we were unable to guess an expression for $C_r$ in terms of $r$ from the five data points). 

Speaking of the OEIS, currently only the cases $r=1$ (of course!) and $r=2$ ([OEIS], sequence A220097) are there.  We hope to
remedy this soon, and enter, at least, $w_r(n)$ for $3 \leq r \leq 5$. 
For $r=3$, although it they are not in OEIS, the first ten terms are already in cyberspace (more precisely, in Lara Pudwell's website).

{
\bf  Some Crucial Background and Zeilberger's  Beautiful Snappy Proof that 123-Avoiding Words are Equinumerous with 132-Avoiding Words
}

Burstein [B] proved that the number of {\it all} words in a given (ordered) alphabet of a given length $n$
avoiding $123$ is the same as the number of words avoiding $132$, and hence, via trivial symmetry, all
patterns of length $3$ have the same enumeration. The stronger result that this is still true if one
specifies the number of occurrences of each letter was 
first proved in [AAAHH], but the {\it proof from the book} appeared in the half-page {\it gem}, [Z3].
Since this lovely proof deserves to be better known, we reproduce it here.

{\it
Define a mapping $F$ on a word $w$ in the alphabet $\{1,2, \dots, n\}$
recursively as follows. If $w$ is empty, then $F(w):=w$.
Otherwise, $i:=w_1$, and let $W$ be the word obtained from
$w$ by first beheading it, and then replacing all letters
larger than $i+1$ by $i+1$, and let $s$ be the sub-sequence
of $w$ obtained by deleting the letters $\leq i$.
Let $\bar{s}$ be the reverse of $s$.
Let $V:=F(W)$, 
and let $U$ be the word obtained from $V$
by replacing (in order) the letters that are $i+1$ by the
members of $\bar{s}$. Finally let $F(w):=iU$.

$F$ is an involution that sends 123-avoiding words to 132-avoiding
ones, and vice versa. This follows from the fact that
$s$ above is non-increasing and non-decreasing respectively.
Hence, for any vector of non-negative integers
$(a_1, \dots, a_n)$ 
amongst the $(a_1+ \dots +a_n)!/(a_1! \cdots a_n!)$ words
with $a_1$ $1$'s, \dots, $a_n$ $n$'s,
the number of those that avoid the pattern $123$ equals the
number of those that avoid $132$,

It also follows that we have a quick recurrence 
that enables us to compute the number
of such words, which we will call $A(a_1, \dots, a_n)$:
$$
A(a_1, \dots, a_n)=\sum_{i=1}^{n} A(a_1,\dots, a_{i-1},a_i -1, a_{i+1}+ \dots + a_n ) \quad.
$$
}

Another important consequence (which also follows from the Robinson-Schenstead-Knuth algorithm) is that
$A(a_1, \dots, a_n)$ is {\it symmetric} in its arguments.

Because of the  equinumeracy of all patterns of length $3$, we can consider $231$-avoiding words.

{\bf Important Definitions}

Let $\W_r(n)$ be the set of $231$-avoiding words in the alphabet $\{1, \dots, n\}$
with exactly $r$ occurrences of each letter.

Also, let $w_r(n)$ be the number of elements of $\W_r(n)$.

Define the `global set'
$$
\W_r \, := \, \bigcup_{n=0}^{\infty} \W_r (n) \quad .
$$
Let $g_r(x)$ be its {\it weight enumerator}  with respect to the weight $w \rightarrow x^{length(w)}$.
Note that $g_r(x)=f_r(x^r)$, where $f_r(x)$ is the generating function of the sequence $w_r(n)$,
$$
f_r(x):=\sum_{n=0}^{\infty} w_r(n) x^n \quad .
$$

We will soon show how, for any {\it specific}, given, positive integer $r$,
to obtain an algebraic equation (i.e. a polynomial $P_r(x,F)$ with integer coefficients such that
$P_r(x,f_r(x))=0$), but let's srart with some {\it warm-ups}.

{\bf First Warm-Up: r=1}

$\W_1$ is the set of {\it all} permutations (of any length!) that avoid the
pattern $231$. Let the weight of a permutation $\pi$ be $x^{length(\pi)}$.
Consider any member $\pi$ of that set. It may happen to be the empty permutation, of course (weight $1$), or else it has a largest
element; let's call that element $n$. All the entries to the left of $n$ must be {\bf smaller} than
all the elements to the right of $n$ (or else a $231$ pattern would emerge),
and each portion must be $231$-avoiding in its own right.
If the location of $n$ is at the $i$-th place, then the portion to the left of $n$ is 
a $231$-avoiding permutation  of $\{1, \dots, i-1 \}$ and the portion to the right is a
$231$-avoiding permutation  of $\{i, \dots, n-1 \}$. Conversely, if $\pi_1$ and $\pi_2$ are
$231$-avoiding permutations of $\{1, \dots, i-1 \}$ and  $\{i, \dots, n-1 \}$ respectively,
then $\pi_1 n \pi_2$ is a $231$-avoiding permutation of length $n$, since no trouble can arise
by joining them. Hence, 
$$
f_1(x)=1+xf_1(x)^2 \quad,
$$
giving the good-old Catalan numbers.

{\bf Second Warm-Up: r=2}

The following argument is inspired by the beautiful proof in [CDZ], but is phrased in such a way
that will make it transparent how to generalize it for general $r$.

Let $g(x)$ be the weight-enumerator  of $\W_2$ . Recall that $\W_2$ is the set of all 
$231$-avoiding words whose letters consist of
$\{ 1,1, \dots, n,n \}$  for some $n \geq 0$, and the weight is $x^{length(w)}=x^{2n}$).

(Note that $g(x)=f_2(x^2))$, so once we have $g(x)$ we will have $f_2(x)$ immediately.)

Consider a typical member of $\W_2$, and let $n$ be its largest element
(i.e. it is of length $2n$). Let $i$ be the location of the {\bf leftmost} occurrence of
$n$. Notice,  just as before, that the entries to the left of that first $n$ must be $\leq$ the entries
to the right of that $n$, and each portion is $231$-avoiding in its own right, and conversely, if
you place such $231$-avoiding words with these entries to the left and right of that leftmost $n$,
you will not cause any trouble, and get a $231$-avoiding word whose entries are $\{1,1,2,2, \dots, n,n\}$.

{\bf Case I}: $i$ is odd, i.e. $i=2j+1$. 

Then the entries to the left of that first $n$ are $\{1,1, \dots, j,j \}$ and the entries to the right are
$\{j+1,j+1, \dots, n-1,n-1, n \}$. The generating function of the left part is our $g(x)$, but the
entries to the right are a new combinatorial creature: a $231$-avoiding word with all the letters
occurring twice, except for one of them (which by symmetry can be taken to be `$1$') that only occurs once. So let's give the set $\W_2$ the new name $\W^{(0,0)}_2$, and let
$\W^{(1,0)}_2$  be the union of the sets of $231$-avoiding words on $\{1,2,2,3,3,\dots, n,n\}$, for all $n \geq 0$.
Let $g^{(1,0)}(x)$  be  its weight-enumerator. Hence the total weight-enumerator of Case I is

$$
x \, g^{(0,0)}(x) \, g^{(1,0)}(x)  \quad .
$$

(The $x$ in front corresponds to the first $n$ separating the two parts).

We will deal with  $g^{(1,0)}_2(x)$ in due course, but now let's proceed to Case II.

{\bf Case II}: $i$ is even, i.e. $i=2j$. 

Once again let its length be $2n$ (so the largest entry is $n$).
The entries to the left of that first $n$ are $\{1,1, \dots, j-1,j-1,j\}$, and the entries to the right are
$\{j, j+1,j+1, \dots, n\}$. The generating function of the left part is the already familiar $g^{(1,0)}(x)$, but the
right part is a new combinatorial creature; namely, a $231$-avoiding word with all the letters
occurring twice, {\bf except} for {\it two} of them (that by symmetry may be taken to be the smallest and the largest)
that only occur {\it once}. Let's call this set
$\W_2^{(1,1)}$,  and its weight-enumerator $g^{(1,1)}(x)$.
Hence the total weight of Case II is $x \, g^{(1,0)}(x) \, g^{(1,1)}(x)$.

Combining the two cases, plus the empty permutation, leads to the following equation
$$
g^{(0,0)}(x)=1 \, + \, x \, g^{(0,0)}(x) \, g^{(1,0)}(x) \, + \, x g^{(1,0)}(x) \, g^{(1,1)}(x) \quad .
\eqno(Eq00)
$$

We have two new {\it uninvited guests},  $g^{(1,0)}(x)$ and $g^{(1,1)}(x)$. Using the
same reasoning as above, the readers are welcome to convince themselves that
$$
g^{(1,0)}(x) \, = \, x g^{(0,0)}(x)^2 \, +  \, x g^{(1,0)}(x)^2 \quad ,
\eqno(Eq10)
$$
$$
g^{(1,1)}(x) \, =\, x  \, g^{(0,0)}(x) \,  g^{(1,0)}(x) \, + \, x \, g^{(1,0)}(x) \, (1+ g^{(1,1)}(x)) \quad .
\eqno(Eq11)
$$

Solving this {\it algebraic scheme}, a system of three algebraic equations $\{Eq00,Eq10,Eq11\}$ in the
three `{\it unknowns}' $\{g^{(0,0)}(x) \,, \, g^{(1,0)}(x) \, , \, g^{(1,1)}(x) \}$, using Gr\"obner bases
(or, in this simple case it could be easily done by hand) gives an algebraic equation satisfied by
$g^{(0,0)}(x)$, and hence, after replacing $x^2$ by $x$, the [CDZ] equation for $f_2(x)$ mentioned above:
$$
1- \left( 2\,x + 1 \right) {f_2(x)}^{2}+x \left( x+4 \right) {f_2(x)}^{4}=0 \quad .
$$

{\bf The General Case}

For $0 \leq i \leq j \leq r-1$ and $n \geq 0$, let $\W_r^{(i,j)}(n)$ be the set of $231$-avoiding words of
length $rn+i+j$, in the alphabet $\{1,2, \dots, n,n+1,n+2\}$, 
with $i$ occurrences of the letter `$1$',  $j$ occurrences of `$n+2$',
and exactly $r$ occurrences of  the other $n$ letters (i.e. $2,3, \dots, n+1$), 
and let  $\W_r^{(i,j)}$ be the union of  $\W_r^{(i,j)}(n)$ over all $n \geq 0$.

By symmetry they have the same weight-enumerator
if {\it any} two letters have $i$ and $j$ occurrences  respectively, and the remaining letters each occur exactly $r$ times.

Using the same logic as above, we have the following ${{r+1} \choose {2}}$ equations,
for  $0 \leq i \leq j \leq r-1$, where below we make the convention that
if $r>s$ then $g^{(r,s)}=g^{(s,r)}$.
$$
g^{(i,j)}(x)=
\delta_{i,0}\delta_{j,0} \, + \,
x \, \sum_{t=0}^{r-1} g^{(i,t)}(x) \, g^{( (r-t) \,mod \, r \, , \, (j-1) \, mod \, r)}(x) \, + \,
\sum_{m=0}^{i-1} x^{m+1} \, g^{(i-m \, , \, j-1)}(x) \quad .
$$

By {\it eliminating} $g^{(0,0)}(x)$, and replacing $x^r$ by $x$, we get the equation of our object of desire $f_r(x)$.
In fact, this equation would have several solutions, and  the right one is picked by plugging in the first few terms.

{\bf Guessing Linear Recurrences for our sequences}

Now that we know, even without WZ-theory, that for every positive integer $r$, the generating function $f_r(x)$ is $D$-finite, since
it has the much stronger property of being algebraic, we immediately know that the sequence itself, $\{ w_r(n) \}$,
is $P$-recursive in the sense of Stanley[S]; in other words, it satisfies {\it some}
homogeneous linear recurrence equation with {\it polynomial} coefficients.

With a very large computer, one should be able to get the algebraic equation for quite a few $r$, and then
use Comtet's algorithm (built-in in the Maple package {\tt gfun}, procedure  {\tt algeqtodiffeq}
followed by  procedure {\tt diffeqtorec}), to get a rigorously derived recurrence. Alas, because
our system has $(r+1)r/2$ algebraic equations, and Gr\"obner bases are notoriously slow, we were only
able to do two new cases explicitly, namely $r=3$ and $r=4$,  mentioned above. But now that we know {\it for sure}
that such recurrences exist, and it is easy to find a priori bounds for the order, it is easy
to justify these empirically-derived recurrences,  a posteriori.

But in order to guess complicated linear recurrences, one needs lots of data. Our algebraic scheme
implies very fast {\it nonlinear} recurrences for the coefficients of $g^{(i,j)}(x)$, and in particular
for  $g^{(0,0)}(x)$, our primary interest. These turn out to be much faster than the `vanilla' linear recurrence
for $A(a_1, \dots, a_n)$ mentioned above.

{\bf The Maple package Words123}

Everything (and more!) is implemented in the Maple package {\tt Words123}, available directly from

{\tt http://www.math.rutgers.edu/\~{}zeilberg/tokhniot/Words123} \quad ,

or via the front of this article

{\tt http://www.math.rutgers.edu/\~{}zeilberg/mamarim/mamarimhtml/words123.html} \quad ,

that also contains some sample input and output files.

{\bf The recurrences for $1 \leq r \leq 3$}

For $r=1$ we get the good-old Catalan numbers
$$
-2\,{\frac { \left( 1+2\,n \right) w_{{1}} \left( n \right) }{n+2}}+w_{{1}} \left( n+1 \right) =0 \quad .
$$

For $r=2$ we get a new proof of the [GGHP] conjecture (first proved in [CDZ])
$$
-3\,{\frac { \left( 7\,n+12 \right)  \left( 1+2\,n \right)  \left( 1+n \right) w_{{2}} \left( n \right) }{ \left( 2\,n+5 \right)  \left( 7\,n+5 \right) 
 \left( n+2 \right) }}
- {\frac { \left( 528+1426\,n+1215\,{n}^{2}+329\,{n}^{3} \right) w_{{2}} \left( n+1 \right) }{ 2 \left( 2\,n+5 \right)  \left( 7\,n+
5 \right)  \left( n+2 \right) }}+w_{{2}} \left( n+2 \right) =0 \quad .
$$
\vfill\eject
For $r=3$ we get
$$
-{\frac {64}{3}}\,{\frac { \left( 4\,n+1 \right)  \left( 2\,n+3 \right)  \left( 4\,n+3 \right)  \left( 14\,n+25 \right)  \left( n+1 \right) w_{{3}} \left( n
 \right) }{ \left( 3\,n+5 \right)  \left( 1+2\,n \right)  \left( 3\,n+7 \right)  \left( 14\,n+11 \right)  \left( n+2 \right) }}
$$
$$
-\frac{8}{3} \cdot{\frac { \left( 3975+
20322\,n+39676\,{n}^{2}+37144\,{n}^{3}+16736\,{n}^{4}+2912\,{n}^{5} \right) w_{{3}} \left( n+1 \right) }{ \left( 3\,n+5 \right)  \left( 1+2\,n \right) 
 \left( 3\,n+7 \right)  \left( 14\,n+11 \right)  \left( n+2 \right) }}+w_{{3}} \left( n+2 \right) =0
\quad .
$$

See the output file

{\tt http://www.math.rutgers.edu/\~{}zeilberg/tokhniot/oWords123c} \quad

for the recurrences for $w_4(n)$ and $w_5(n)$.

{\bf The Asymptotics for $1 \leq r \leq 5$}

$$
w_1(n) \, = \,
{\frac {1}{\sqrt {\pi }}} \, \cdot \, {4}^{n} \,
\, \cdot \, n^{-\frac{3}{2}} \,
\left( 1-{\frac {9}{8}}\,{n}^{-1}+{\frac {145}{128}}\,{n}^{-2}-{\frac {1155}{1024}}\,{n}^{-3} 
\, + \, O( n^{-4}) \,
\right)
\quad , 
$$
$$
w_2(n) \, = \,
{\frac {1}{\sqrt {\pi }}} \, \cdot \,
{\frac {3 \, \sqrt {3} }{7 \sqrt{7}}} \, \cdot \, {12}^{n} 
\, \cdot \, n^{-\frac{3}{2}} \,
\left( 1-{\frac {249}{392}}\,{n}^{-1}+{\frac {13255}{43904}}\,{n}^{-2}-{\frac {
2674485}{17210368}}\,{n}^{-3} 
\, + \, O( n^{-4} ) \,
\right)
\quad , 
$$
$$
w_3(n) \, = \,
\frac{1}{\sqrt {\pi }} \, \cdot \, \frac{1}{8} \, \cdot  \,
{32}^{n} \, 
\, \cdot \, n^{-\frac{3}{2}} \,
\left( 1-{\frac {33}{64}}\,{n}^{-1}+{\frac {1105}{8192}}\,{n}^{-2}-{\frac {27195}{524288}}\,{n}^{-3} 
\, + \, O( n^{-4} \right) 
\quad ,
$$
$$
w_4(n) \, =  \frac{1}{\sqrt{\pi}} \, \cdot \, \frac{1}{6 \sqrt{6} } \, \cdot \,
{80}^{n}  
\, \cdot \, n^{-\frac{3}{2}} \,
\left( 1-{\frac {23}{48}}\,{n}^{-1}+{\frac {1621}{23040}}\,{n}^{-2}-{\frac {339199}{16588800
}}\,{n}^{-3} 
\, + \, O( n^{-4}) \,
\right ) \quad ,
$$
$$
w_5(n) \, = \,
\frac{1}{\sqrt {\pi }} \, \cdot \,
\frac{3 \sqrt {3}}{125} \, \cdot {192}^{n}  \,
\, \cdot \, n^{-\frac{3}{2}} \,
\left( 1-{\frac {471}{1000}}\,{n}^{-1}+{\frac {389141}{10000000}}\,{n}^{-2}-{\frac {
162387477}{50000000000}}\,{n}^{-3} 
\, + \, O( n^{-4}) \,
\right )  \quad .
$$
\bigskip
{\eightrm Warning: the above asymptotic expressions are fully rigorous except for the constants in front, which are only conjectured}.

{\bf References}

[AAAHH] M. H. Albert, M. Aldred, M. D. Atkinson, C. C. Handley and D. A. Holton,
{\it Permutations of a multiset avoiding permutations of length $3$}, Europ. J. Comb. {\bf 22} (2001), 1021-1031. \hfill \break
{\tt http://reflect.otago.ac.nz/staffpriv/mike/Papers/Multiperms/Multiperms.pdf} \quad .

[B] A. Burstein, {\it ``Enumeration of words with forbidden patterns''}, PhD Thesis, University of Pennsylvania, 1998.

[CDZ]   W. Y. C. Chen, A. Y. L. Dai and R. D. P. Zhou, {\it Ordered Partitions Avoiding a Permutation of Length 3}, 
to appear in European J. of Combinatorics.
\hfill\break
{\tt http://arxiv.org/abs/1304.3187} \quad .

[Co] L. Comtet, {\it Calcul pratique des coefficients de Taylor d'une fonction alg\'ebrique}, Enseignement Mathematique
{\bf 10} (1964), 267-270.

[GGHP] A. Godbole, A.Goyt, J. Herdan, and L. Pudwell, {\it Pattern Avoidance in Ordered Set Partitions}, 
Ann. Comb. {\bf 18} (2014), 429-445. \hfill\break
{\tt http://arxiv.org/abs/1212.2530} \quad .

[HM] S. Heubach and T. Mansour, ``{\it Combinatorics of compositions and words}'', 
Discrete Mathematics and its Applications (Boca Raton), CRC Press, Boca Raton, FL, 2010 \quad .

[Ka] Anisse Kasraoui, {\it Pattern avoidance in ordered set partitions and words}, Advances in Applied Mathematics
Oct. 2014.
\hfill\break
{\tt http://arxiv.org/abs/1307.0495} .

[Ki] S. Kitaev, ``{\it Patterns in permutations and words}'', Springer Verlag (EATCS monographs in Theoretical Computer Science book series), 2011. 

[OEIS] The OEIS Foundation, Sequence A220997, {\tt https://oeis.org/A220097} \quad .

[P] L. Pudwell, 
{\it Enumeration schemes for words avoiding permutations}. 
In ``Permutation Patterns'' ({\bf 2010}), S. Linton, N. Ruskuc, and V. Vatter, Eds., 
vol. {\bf 376} of London Mathematical Society Lecture Note Series, Cambridge University Press, pp. 193-211. Cambridge: Cambridge University Press.

[S] R. Stanley, {\it Differentiably finite power series}, European J. Combinatorics {\bf 1} (1980), 175-188. \hfill\break
{\tt http://www-math.mit.edu/\~{}rstan/pubs/pubfiles/45.pdf} \quad .

[Wiki] Wikipedia contributors (most notably Vince Vatter),
{\it `` Permutation pattern''}, Wikipedia, The Free Encyclopedia, 19 Sep. 2014. Web. 7 Nov. 2014. \hfill\break
{\tt https://en.wikipedia.org/wiki/Permutation\_pattern} \quad .

[Z1] D. Zeilberger, {\it A Holonomic Systems Approach To Special Functions},
J. Computational and Applied Math {\bf 32} (1990), 321-368. \hfill\break
preprint version: {\tt http://www.math.rutgers.edu/\~{}zeilberg/mamarim/mamarimPDF/holonomic.pdf}

[Z2] D. Zeilberger,
{\tt AsyRec: A Maple package for Computing the Asymptotics of Solutions of Linear Recurrence Equations with Polynomial Coefficients},
Personal Journal of Shalosh B. Ekhad and Doron Zeilberger, April 4, 2008. \hfill\break
{\tt http://www.math.rutgers.edu/\~{}zeilberg/mamarim/mamarimhtml/asy.html}

[Z3] D. Zeilberger,
{\it A Snappy Proof That 123-Avoiding Words are Equinumerous With 132-Avoiding Words},
Personal Journal of Shalosh B. Ekhad and Doron Zeilberger, April 11, 2005. \hfill\break
{\tt http://www.math.rutgers.edu/\~{}zeilberg/mamarim/mamarimhtml/a123.html} \quad .

\vfill\eject

[Z4] D. Zeilberger,
{\it On Vince Vatter's Brilliant Extension of Doron Zeilberger's Enumeration Schemes for Counting Herb Wilf's Classes},
Personal Journal of Shalosh B. Ekhad and Doron Zeilberger, Dec. 29, 2006. \hfill\break
{\tt http://www.math.rutgers.edu/\~{}zeilberg/mamarim/mamarimhtml/vatter.html} \quad .

\bigskip
\hrule
\bigskip
Nathaniel Shar, Department of Mathematics, Rutgers University (New Brunswick), Hill Center-Busch Campus, 110 Frelinghuysen
Rd., Piscataway, NJ 08854-8019, USA. \hfill\break
Email: {\tt nshar at math dot rutgers dot edu}   \quad .
\bigskip
Doron Zeilberger, Department of Mathematics, Rutgers University (New Brunswick), Hill Center-Busch Campus, 110 Frelinghuysen
Rd., Piscataway, NJ 08854-8019, USA. \hfill\break
Email: {\tt zeilberg at math dot rutgers dot edu}   \quad .
\bigskip
\hrule
\bigskip
Nov. 18, 2014

\end